\documentclass[oneside, reqno]{amsart}

\usepackage{xypic}
\input xy
\xyoption{all}
\usepackage{epsfig}
\usepackage{amsthm}
\usepackage{amssymb}
\usepackage{amsmath}
\usepackage{amscd}
\usepackage{color}
\usepackage[T1]{fontenc}
\usepackage{stackengine}
\usepackage{upgreek}
\usepackage[font=scriptsize]{caption}
\usepackage{wrapfig}
\usepackage[utf8x]{inputenc}
\stackMath

\newcommand{\bg}{\begin{equation}}
\newcommand{\ed}{\end{equation}}
\newcommand{\bga}{\begin{eqnarray}}
\newcommand{\eda}{\end{eqnarray}}
\newcommand{\pf}{\textbf{Proof:\ }}

\def\cbdu{\par{\raggedleft$\Box$\par}}

\newtheorem {Theorem}  {Theorem}

\numberwithin{Theorem}{section}

\newtheorem {Lemma}[Theorem]  {Lemma}

\theoremstyle{definition}
\newtheorem{Definition}[Theorem]{Definition}
\theoremstyle{remark}
\newtheorem{Remark}[Theorem]{\bf Remark}

\expandafter\chardef\csname pre amssym.def
at\endcsname=\the\catcode`\@ \catcode`\@=11
\def\undefine#1{\let#1\undefined}
\def\newsymbol#1#2#3#4#5{\let\next@\relax
 \ifnum#2=\@ne\let\next@\msafam@\else
 \ifnum#2=\tw@\let\next@\msbfam@\fi\fi
 \mathchardef#1="#3\next@#4#5}
\def\mathhexbox@#1#2#3{\relax
 \ifmmode\mathpalette{}{\m@th\mathchar"#1#2#3}%
 \else\leavevmode\hbox{$\m@th\mathchar"#1#2#3$}\fi}
\def\hexnumber@#1{\ifcase#1 0\or 1\or 2\or 3\or 4\or 5\or 6\or 7\or 8\or
 9\or A\or B\or C\or D\or E\or F\fi}

\font\teneufm=eufm10 \font\seveneufm=eufm7 \font\fiveeufm=eufm5
\newfam\eufmfam
\textfont\eufmfam=\teneufm \scriptfont\eufmfam=\seveneufm
\scriptscriptfont\eufmfam=\fiveeufm

\catcode`\@=\csname pre amssym.def at\endcsname

\newcounter{remark}
\def  \12  {{\frac{1}{2}}}

\def\build#1_#2^#3{\mathrel{\mathop{\kern 0pt#1}\limits_{#2}^{#3}}}

\numberwithin{equation}{section}

\begin{document}

\title[Dyadic MHD model]{Dyadic models for ideal MHD}


\author [Mimi Dai]{Mimi Dai}

\address{Department of Mathematics, Statistics and Computer Science, University of Illinois at Chicago, Chicago, IL 60607, USA}
\email{mdai@uic.edu} 

\author [Susan Friedlander]{Susan Friedlander}

\address{Department of Mathematics, University of Southern California, Los Angeles, CA 90089, USA}
\email{susanfri@usc.edu}


\begin{abstract}

We study two dyadic models for incompressible ideal magnetohydrodynamics, one with a uni-directional energy cascade and the other one with both forward and backward energy cascades. Global existence of weak solutions and local well-posedness are established for both models. In addition, solutions to the model with uni-directional energy cascade associated with positive initial data are shown to develop blow-up at a finite time. Moreover, a set of fixed points is found for each model. Linear instability about some particular fixed points is proved.

\bigskip


\hspace{0.02cm}CLASSIFICATION CODE: 35Q35, 76B03, 76W05.
\end{abstract}

\maketitle

\section{Introduction}

\subsection{Dyadic Euler model}
Dyadic models for hydrodynamics have a long history traced back to the first such models introduced by Desnyanskiy and Novikov \cite{DN}, Gledzer \cite{Gle}, and Ohkitani and Yamada \cite{OY}, usually referred as DN model and GOY model respectively. 
Following these early papers, various models \cite{Bif, CLT, DS, Ob} were developed by simplifying the nonlinearity in different ways.
 Among them, the Katz-Pavl\'ovic (KP) model \cite{KP} was derived applying Littlewood-Paley techniques to the Euler equation
 \begin{equation}\notag
 u_t+(u\cdot\nabla) u+\nabla P=0.
 \end{equation}
 The KP model with external forcing is the following nonlinearly coupled infinite system of ordinary differential equations (ODEs)
\begin{equation}\label{sys-nse}
\frac{d}{dt}a_j +\lambda_j^{\theta}a_ja_{j+1}-\lambda_{j-1}^{\theta}a_{j-1}^2=f_j,
\end{equation}
for $j\geq 0$ and $a_{-1}=0$, 
with $f=(f_0, f_1, f_2, ...)$. The quantity $\frac12a_j^2$ stands for the energy of the $j$-th dyadic shell with wavenumber $\lambda_j=\lambda^j$ for a constant $\lambda>1$. The parameter $\theta$ takes the form $\theta=\frac{5-\delta}{2}$ with $\delta$ being the intermittency dimension for the turbulent fluid (cf. \cite{CD-Kol}). For 3D fluids, $\delta\in[0,3]$ and hence $\theta\in[1, \frac52]$. Evidently smaller $\delta$ corresponds to larger $\theta$. Conceptually a turbulent flow with smaller intermittency dimension $\delta$ is more singular and hence exhibits stronger nonlinear effects.

The main features of (\ref{sys-nse}) include: 
(i) it formally conserves the energy when $f=0$;  (ii) it only takes into account the nonlinear interactions of nearest neighboring shells through a Littlewood-Paley decomposition treatment; (iii) it is equipped with a forward energy cascade mechanism; (iv) spatial structure is oversimplified. Although model (\ref{sys-nse}) has a seemingly simple form, analyzing it is highly non-trivial. Fortunately, two important properties play a vital role in the analysis: positivity of solutions starting from positive data and monotonicity of the rescaled quantity $\lambda_j^{\frac13\theta} a_j$.

There is an extensive literature for the study of model (\ref{sys-nse}),
cf. \cite{CGHPSV} for a thorough description. We highlight some results from the literature on the Euler model without  the intention of providing a complete list. Finite time blow-up was first shown in \cite{KP} and later sharpened in \cite{KZ}. For a closely related vector model, the authors of \cite{FP} also discovered finite time blow-up. A model of a combination of DN type and Obukhov type \cite{Ob} was studied in \cite{JL} and shown to exhibit finite time blow-up when the Obukhov type nonlinearity is sufficiently weak. 
Anomalous energy dissipation and self-similar solutions for the unforced model were studied in \cite{BFM}. The problems of uniqueness, well-posedness and regularity for the inviscid model were addressed in \cite{BFM-uni, BM}. The authors of \cite{CZ} proved regularizing properties of the nonlinear term of (\ref{sys-nse}) due to the forward energy cascade mechanism. The models with stochastic forcing were also investigated in many articles, for instance, see \cite{BFM-stoc1, BFM-stoc2, FGHV, Rom}. 

The results of Cheskidov, Friedlander and Pavl\'ovic given in \cite{CFP1, CFP2} are relevant to our current paper. In \cite{CFP1}, the authors studied (\ref{sys-nse}) with $\theta=1$ and positive forcing. They showed the existence of a unique positive fixed point for the dyadic Euler model (\ref{sys-nse}) with positive initial data. They also established linear stability of the fixed point in the sense that there is no positive eigenvalue for the linearized system about the fixed point and that there exist solutions to the linearized system in $H^s$ with $s<\frac13$ which decay exponentially fast in time. Moreover, they proved that every solution with bounded $H^{1/3}$ norm approaches the fixed point in the energy space $l^2$. Hence solutions with positive initial $l^2$ data blow up in finite time in $H^{1/3}$ norm. In \cite{CFP2} the authors further examined the properties of the unique positive fixed point and showed that the fixed point is an exponential global attractor. Via this property the authors provided another proof of finite time blowup of solutions in the norm $H^s$ with $s\geq\frac13$. We note that persistence of positivity of the solutions plays a crucial role in the analysis of \cite{CFP1, CFP2}. 

\subsection{Dyadic models for magnetohydrodynamics}
The equations for incompressible ideal (inviscid and non-resistive) magnetohydrodynamics (MHD) are the set of partial differential equations
\begin{subequations}
\begin{align}
u_t+(u\cdot\nabla) u-(B\cdot\nabla) B+\nabla P=&\ 0, \label{mhda}\\
B_t+(u\cdot\nabla) B-(B\cdot\nabla) u =&\ 0, \label{mhdb}\\
\nabla \cdot u= 0, \ \ \nabla \cdot B=&\ 0. \label{mhdc}
\end{align}
\end{subequations}
System (\ref{mhda})-(\ref{mhdc}) describes electrically conducting fluids in geophysics and astrophysics when the underlying length scales are very large, and hence the kinetic viscosity and magnetic diffusivity are negligible.  
In (\ref{mhda})-(\ref{mhdc}) the unknown vector fields $u$ and $B$ and the scalar function $P$ denote the fluid velocity, magnetic field and pressure. 

Inspired by results described above for the dyadic model (\ref{sys-nse}) of the Euler equation, we study two particular dyadic models (suggested in \cite{Dai-20}) associated with the PDEs (\ref{mhda})-(\ref{mhdc}) of MHD. 
One such model has a uni-directional energy cascade. The other model has both forward and backward energy cascade mechanisms. The model with uni-directional energy cascade is presented as 
\begin{equation}\label{sys-2}
\begin{split}
\frac{d}{dt}a_j=& -\left(\lambda_j^{\theta}a_ja_{j+1}-\lambda_{j-1}^{\theta}a_{j-1}^2\right)
-\left(\lambda_{j}^{\theta}b_jb_{j+1}-\lambda_{j-1}^{\theta}b_{j-1}^2\right) +f_j,\\
\frac{d}{dt}b_j= &\ \lambda_j^{\theta}a_jb_{j+1}-\lambda_{j}^{\theta}b_ja_{j+1},
\end{split}
\end{equation}
for $j\geq0$ with $a_{-1}=b_{-1}=0$ and $f=(f_0, f_1, f_2, ...)$. For positive solutions, the forward energy cascade within the dynamics can be illustrated below
\[ \cdot\cdot\cdot  \longrightarrow  \ \ a_{j-1} \ \ \longrightarrow \ \ a_j \ \ \longrightarrow \ \ a_{j+1} \ \ \longrightarrow \ \cdot\cdot\cdot\]
\[\hspace{-4mm}\downarrow \ \ \ \ \ \ \nearrow \ \  \ \ \downarrow \ \ \hspace{4mm} \nearrow \ \  \ \ \downarrow\]
\[ \cdot\cdot\cdot  \longrightarrow  \ \ b_{j-1} \ \ \longrightarrow \ \ b_j \ \ \longrightarrow \ \ b_{j+1} \ \ \longrightarrow \ \cdot\cdot\cdot\]
The model with both forward and backward energy cascades is 
\begin{equation}\label{sys-1}
\begin{split}
\frac{d}{dt}a_j=&-\left(\lambda_j^{\theta}a_ja_{j+1}-\lambda_{j-1}^{\theta}a_{j-1}^2\right)+\left(\lambda_{j}^{\theta}b_jb_{j+1}-\lambda_{j-1}^{\theta}b_{j-1}^2\right) +f_j,\\
\frac{d}{dt}b_j= & -\left(\lambda_j^{\theta}a_jb_{j+1}-\lambda_{j}^{\theta}b_ja_{j+1}\right),
\end{split}
\end{equation}
with the bi-directional energy transfer illustrated as, 
\[ \cdot\cdot\cdot  \longrightarrow  \ \ a_{j-1} \ \ \longrightarrow \ \ a_j \ \ \longrightarrow \ \ a_{j+1} \ \ \longrightarrow \ \cdot\cdot\cdot\]
\[\hspace{-4mm}\uparrow \ \ \ \ \ \ \swarrow \ \  \ \ \uparrow \ \ \hspace{4mm} \swarrow \ \  \ \ \uparrow\]
\[ \cdot\cdot\cdot  \longrightarrow  \ \ b_{j-1} \ \ \longrightarrow \ \ b_j \ \ \longrightarrow \ \ b_{j+1} \ \ \longrightarrow \ \cdot\cdot\cdot\]
The quantities $\frac12a_j^2$ and $\frac12b_j^2$ in (\ref{sys-2}) and (\ref{sys-1}) represent the kinetic energy and magnetic energy in the $j$-th shell, respectively. One can see that only interactions with the nearest neighbor shells are taken into account in the modeling. Both models (\ref{sys-2}) and (\ref{sys-1}) preserve the most essential feature of the original dynamics 
of (\ref{mhda})-(\ref{mhdc}), namely the total energy 
\begin{equation}\notag
E(t)=\frac12\sum_{j\geq 1}\left(a_j^2(t)+b_j^2(t)\right)
\end{equation}
is conserved formally. In addition, the cross helicity defined by
\[H^c(t)=\sum_{j=0}^\infty a_jb_j\]
is invariant for strong solutions of the model (\ref{sys-1}), but not for model (\ref{sys-2}).

It is obvious that  eliminating $b$ (i.e. $b_j=0$ for $j\geq0$) reduces systems (\ref{sys-2}) and (\ref{sys-1}) to the fluid model (\ref{sys-nse}). Beside sharing all of the similar difficulties arising from the nonlinearity in (\ref{sys-nse}), the MHD models (\ref{sys-2}) and (\ref{sys-1}) exhibit extra obstacles due to the interactions between the fluid and magnetic fields. In particular, it is unlikely that solutions with positive initial data stay positive for all the time in a general setting. It is also not clear whether the quantities $\lambda_j^{\frac13\theta} a_j$ and $\lambda_j^{\frac13\theta} b_j$ exhibit certain monotonicity. On the other hand, the conservation of positivity is not a natural property, especially in the context of MHD models. Numerical simulations suggest solutions with positive initial data can become negative.

We take a brief look at the energy transfer mechanisms and get some basic sense of the difficulties.  
Denote the flux by $\Pi_j=\lambda_{j}^\theta \left(a_{j}^2+b_{j}^2\right)a_{j+1}$ for $j\geq 0$. The rate of change of the total energy of the $j$-th shell for (\ref{sys-2}) is given by
\begin{equation}\notag
\begin{split}
\frac12 \frac{d}{dt} \left(a_j^2+b_j^2\right)=&\ \lambda_{j-1}^\theta \left(a_{j-1}^2+b_{j-1}^2\right)a_j-\lambda_j^\theta \left(a_j^2+b_j^2\right)a_{j+1}+f_ja_j\\
=&\ \Pi_{j-1}-\Pi_j+f_ja_j,
\end{split}
\end{equation}
where $\Pi_{j-1}$ is the flux coming from the previous shell and $\Pi_j$ is the flux escaping to the next shell provided $a_j, a_{j+1}\geq0$, and vice versa if $a_j, a_{j+1}\leq 0$. It thus represents a uni-directional energy cascade mechanism. Nevertheless, since the conservation of positivity is not valid in general setting, analyzing the energy transfer is not trivial. Moreover, there are obstacles to understand the energy flux $b_j^2a_{j+1}$ involved with interactions between the velocity and the magnetic fields.
The situation is even more subtle for system (\ref{sys-1}). Denote $\tilde\Pi_j=\lambda_{j}^\theta \left(a_{j}^2-b_{j}^2\right)a_{j+1}$ for $j\geq 0$. Then for (\ref{sys-1}) the energy change rate of the $j$-th shell is
\begin{equation}\notag
\frac12 \frac{d}{dt} \left(a_j^2+b_j^2\right)= \tilde\Pi_{j-1}-\tilde\Pi_j+f_ja_j.
\end{equation}
We note that the sign of the flux $\tilde \Pi_j$ does not only depend on $a_{j+1}$ but also on the size of $|a_j|$ and $|b_j|$. Hence it is challenging to obtain a rather comprehensive picture of the energy transfer for (\ref{sys-1}).

Despite the obstacles described above, the authors \cite{DF} were able to obtain results on uniqueness and non-uniqueness for the models (\ref{sys-2}) and (\ref{sys-1}) with viscous linear terms $\lambda_j^2a_j$ and $\lambda_j^2b_j$. More precisely, uniqueness of Leray-Hopf solution was established for $\theta\leq 2$; while non-unique Leray-Hopf solutions were constructed for $\theta>2$. We remark that positivity is not known to hold for the viscous dyadic MHD models. The techniques used in \cite{DF} are independent of the positivity property.

We note that dyadic models for the MHD turbulence were introduced a few decades ago in the physics community, see \cite{Bis94, Frik, GLPG}. A majority of contributions, for instance in \cite{AMP, AF, ALFS, BSDP}, concern numerical studies of various dyadic models. In fact, one of the motivations to develop dyadic models for the MHD is that direct numerical simulations for the original dynamics face serious computing limitations. More complex dyadic models have been proposed to include non-local interactions \cite{FSS} and anisotropy \cite{NMCV} to understand the MHD turbulence. Although the dyadic models do not preserve geometric features due to the lack of spatial structure, numerics for these models show intermittency statistics and chaotic behaviour, which are in agreement with experimental MHD turbulence.    
For a full background in this regard, the reader may consult the recent review article \cite{PSF} and references therein.

The main objective of the present article is to examine the inviscid systems (\ref{sys-2}) and (\ref{sys-1}). Global weak solutions and local well-posedness are obtained by employing standard arguments. In the context of the dyadic MHD model (\ref{sys-2}), we show finite time blow-up for (\ref{sys-2}) by constructing an appropriate Lyapunov function. We note that it remains an open question whether a solution of (\ref{sys-1}) with bi-directional energy cascade develops blow-up in finite time.

A set of fixed points for (\ref{sys-2}) and (\ref{sys-1}) are found with an explicit form. Specifically, assuming $f$ takes the form
\[f_0>0, \ \ \mbox{and} \ \ f_j=0, \ \ \forall \ \ j\geq 1,\]
we find infinitely many steady states $(\bar a, \bar b)$ of (\ref{sys-2}) in $l^2\times l^2$ which satisfy 
\[\bar a_j^2+\bar b_j^2=\lambda^{\frac13\theta} f_0 \lambda_j^{-\frac23\theta}, \ \ \forall \ \ j\geq 0,\]
and infinitely many steady states $(\bar a, \bar b)$ of (\ref{sys-1}) satisfying 
\[\bar a_j^2-\bar b_j^2=\lambda^{\frac13\theta} f_0 \lambda_j^{-\frac23\theta}, \ \ \forall \ \ j\geq 0.\]
In stark contrast with the dyadic Euler model studied in \cite{CFP1, CFP2}, it is much harder to analyze the properties of the fixed points of systems (\ref{sys-2}) and (\ref{sys-1}). One obvious reason is that the fixed points are not unique. The more sophisticated difficulty comes from the interactions between the velocity and magnetic fields. Our understanding of the energy transfer from the $j$-th shell of velocity to the $j$-th shell of magnetic field ($a_j$ to $b_j$) is limited at the moment.  Although it is challenging to show stability or instability of the fixed points in general setting, we are able to establish linear instability results for some special fixed points. In particular, for the fixed point with zero magnetic field $(\bar a, 0)$, we observe that the linearized system is decoupled between the velocity and magnetic field and we prove that it is linearly unstable. More specifically, if $\bar a$ is positive, the velocity component is stable while the magnetic field component is unstable; if $\bar a$ is negative, both of the velocity and magnetic field are unstable.

\subsection{Organisation of the paper}

\begin{itemize}
\item Section \ref{sec-pre} consists of notations and definitions of solutions for the dyadic systems (\ref{sys-2}) and (\ref{sys-1}).
\item Section \ref{sec-weak} gives a proof of global existence of weak solutions for any initial data with finite energy. This result holds for both systems (\ref{sys-2}) and (\ref{sys-1}).
\item Section \ref{sec-well} gives a proof of local well-posedness in $H^s$ for $s\geq \theta$. This result holds for both systems (\ref{sys-2}) and (\ref{sys-1}).
\item In Section \ref{sec-blow} we show that the system (\ref{sys-2}) starting from positive initial data blows up in finite time in $H^s$ for $s>\frac{1}{3}\theta$.
\item Section \ref{sec-fixed} describes the fixed point sets for systems (\ref{sys-2}) and (\ref{sys-1}).
\item In Section \ref{sec-stab} we examine the linear instability for perturbations about two specific fixed points.
\end{itemize}

\bigskip

\section{Preliminaries}
\label{sec-pre}

\subsection{Notation and definition of solutions}
We will often use $C_0$ and $c, c_0, c_1, ...$ to denote constants which may vary from line to line. They are universal constants unless specified otherwise. We say $a=(a_0, a_1, a_2, ...)$ is positive if all of the components are positive.

Parallel to the energy space $L^2$ and Sobolev space $H^s$ for functions of space and time, we need to refer an analogue of such spaces for sequences by $l^2$ and $H^s$ (the same notation), with $l^2$ endowed with 
 the standard scalar product and norm
\[(u,v):=\sum_{n=1}^\infty u_nv_n, \ \ \ |u|:=\sqrt{(u,u)},\]
and $H^s$ equipped with the scaler product 
\[(u,v)_s:=\sum_{n=1}^\infty \lambda_n^{2s}u_nv_n \ \ \ \ \mbox{with} \ \ \lambda_n=\lambda^n\]
and the norm
\[\|u\|_{H^s}=\|u\|_s:=\sqrt{(u,u)_s}.\]
We define the strong distance $\mathrm d_s$ and weak distance $\mathrm d_w$ on $l^2$ as,
\[\mathrm d_{\mathrm s}(u,v):=|u-v|, \ \ \ \mathrm d_{\mathrm w}(u,v):=\sum_{n=1}^\infty \frac{1}{\lambda^{n^2}}\frac{|u_n-v_n|}{1+|u_n-v_n|}, \ \ u,v\in l^2.\]
Naturally, $\mathrm d_{\mathrm w}$ generates a weak topology on any bounded subset of $l^2$. Weak convergence of a bounded sequence $\{u^k\}\subset l^2$ to $u\in l^2$ is understood in the usual way, 
\[\mathrm d_{\mathrm w}(u^k,u) \to 0 \ \ \ \mbox{as} \ \ \ k\to \infty. \] 
The functional spaces $C([0,T]; l^2_{\mathrm w})$ and $C([0,\infty); l^2_{\mathrm w})$ are defined as  
\[C([0,T]; l^2_{\mathrm w}):= \{u(\cdot): [0,T]\to l^2, \ u_n(t) \ \mbox{is continuous for all} \ n\}\]
endowed with the distance 
\[\mathrm d_{C([0,T]; l^2_{\mathrm w})}(u,v): =\sup_{t\in[0,T]} \mathrm d_{\mathrm w}(u(t), v(t)),\]
and 
\[C([0,\infty); l^2_{\mathrm w}):= \{u(\cdot): [0,\infty)\to l^2, \ u_n(t) \ \mbox{is continuous for all} \ n\}\]
with 
\[\mathrm d_{C([0,\infty); l^2_{\mathrm w})}: =\sum_{T\in \mathbb N}\frac{1}{2^T}\frac{\sup\{\mathrm d_{\mathrm w}(u(t),v(t)): 0\leq t\leq T\}}{1+\sup\{\mathrm d_{\mathrm w}(u(t),v(t)): 0\leq t\leq T\}}.\]

\begin{Definition}\label{def1}
A pair of $l^2$-valued functions $(a(t), b(t))$ defined on $[t_0,\infty)$ is said to be a weak solution of (\ref{sys-2}) or (\ref{sys-1}) if $a_j$ and $b_j$ satisfy (\ref{sys-2}) or (\ref{sys-1})and $a_j, b_j\in C^1([t_0,\infty))$ for all $j\geq0$.
\end{Definition}

\begin{Definition}\label{def2}
A solution $(a(t), b(t))$ of (\ref{sys-2}) or (\ref{sys-1}) is strong on $[T_1, T_2]$ if $\|a\|_{H^1}$ and $\|b\|_{H^1}$ are bounded on $[T_1, T_2]$. A solution is strong on $[T_1, \infty)$ if it is strong on every interval $[T_1, T_2]$ for any $T_2>T_1$.
\end{Definition}

\bigskip

\section{Weak solutions}\label{sec-weak}

In this part, we show the global existence of weak solutions for any data with finite energy by the Galerkin approximating method.

\begin{Theorem}\label{thm-weak}
For any initial data $(a^0, b^0)\in l^2\times l^2$, there exists a weak solution $(a(t), b(t))$ of system (\ref{sys-2}) and (\ref{sys-1}) on $[0,\infty)$ with $a(0)=a^0$ and $b(0)=b^0$.
\end{Theorem}
\pf
We only need to show the existence of a weak solution to system (\ref{sys-2}), since minor modification of the proof works for system (\ref{sys-1}). 
For any fixed integer $k\geq 1$, we consider the truncated variables 
$\{(a^k(t), b^k(t))\}$ given by
\[a^k(t)=\left(a^k_0(t), a^k_1(t), a^k_2(t), ..., a^k_k(t), 0, 0, ...\right), \]
\[b^k(t)=\left(b^k_0(t), b^k_1(t), b^k_2(t), ..., b^k_k(t), 0, 0, ...\right), \]
with $a_j^k(0)=a_j^0$ and $b_j^k(0)=b_j^0$ for $0\leq j\leq k$, and
$\{(a^k(t), b^k(t))\}$ satisfies the finite system
\begin{equation}\label{sys-truncated}
\begin{split}
\frac{d}{dt}a_j^k=&-\lambda_j^{\theta}a_j^ka_{j+1}^k+\lambda_{j-1}^{\theta}(a_{j-1}^k)^2
-\lambda_{j}^{\theta}b_j^kb_{j+1}^k\\
&+\lambda_{j-1}^{\theta}(b_{j-1}^k)^2+f_j,  \ \ \ \ j\leq k-1\\
\frac{d}{dt}b_j^k=&\ \lambda_j^{\theta}a_j^kb_{j+1}^k-\lambda_{j}^{\theta}b_j^ka_{j+1}^k, \ \ \ \  j\leq k-1\\
\frac{d}{dt}a_k^k=&\ \lambda_{k-1}^{\theta}(a_{k-1}^k)^2+\lambda_{k-1}^{\theta}(b_{k-1}^k)^2+f_k,  \\
\frac{d}{dt}b_k^k=&\ 0.  \\
\end{split}
\end{equation}
We notice that the functions on the right hand side of (\ref{sys-truncated}) are continuous in $a^k$ and $b^k$. Hence, there exists a unique solution $(a^k(t), b^k(t))$ to (\ref{sys-truncated}) on $[0,T]$ for arbitrary $T>0$. In order to pass a subsequence of $\{(a^k(t), b^k(t))\}$ to a limit, we will apply the Ascoli-Arzela theorem. It is sufficient to show that the sequence $\{(a^k(t), b^k(t))\}$ is weakly equicontinuous. For any fixed $k\geq 1$, there exists a constant $C_0$ independent of $k$ such that 
\begin{equation}\label{bound}
|a_j^k(t)|\leq C_0, \ \ \ |b_j^k(t)|\leq C_0 \ \ \ \ \forall \ t\in[0,T] \ \ \mbox{and} \ \ \ 0\leq j\leq k.
\end{equation}
It follows from (\ref{sys-truncated}) and (\ref{bound}) that 
\begin{equation}\label{est-ajk}
\begin{split}
&|a_j^k(t)-a_j^k(s)|\\
\leq &\int_s^t\left| -\lambda_j^{\theta}a_j^k(\tau)a_{j+1}^k(\tau)+\lambda_{j-1}^{\theta}(a_{j-1}^k(\tau))^2 \right.\\
& \left.-\lambda_{j}^{\theta}b_j^k(\tau)b_{j+1}^k(\tau) +\lambda_{j-1}^{\theta}(b_{j-1}^k(\tau))^2+f_j\right| \, d\tau\\
\leq & \left(2\lambda_j^\theta C_0^2+2\lambda_{j-1}^\theta C_0^2+f_j\right) |t-s|
\end{split}
\end{equation}
for $0\leq j\leq k$ and $0\leq t\leq s\leq T$, and similarly 
\begin{equation}\label{est-bjk}
\begin{split}
&|b_j^k(t)-b_j^k(s)|\\
\leq &\int_s^t\left| \lambda_j^{\theta}a_j^k(\tau)b_{j+1}^k(\tau)
-\lambda_{j}^{\theta}b_j^k(\tau)a_{j+1}^k(\tau)\right| \, d\tau\\
\leq &\ 2\lambda_j^\theta C_0^2 |t-s|.
\end{split}
\end{equation}
Therefore, we have from (\ref{est-ajk}) and (\ref{est-bjk})
\begin{equation}\notag
\begin{split}
\mathrm d_{\mathrm w} (a^k(t), a^k(s))=&\sum_{j=0}^\infty \frac{1}{\lambda^{j^2}}\frac{|a_j^k(t)-a_j^k(s)|}{1+|a_j^k(t)-a_j^k(s)|}\leq c|t-s|,\\
\mathrm d_{\mathrm w} (b^k(t), b^k(s))=&\sum_{j=0}^\infty \frac{1}{\lambda^{j^2}}\frac{|b_j^k(t)-b_j^k(s)|}{1+|b_j^k(t)-b_j^k(s)|}\leq c|t-s|,
\end{split}
\end{equation}
for an absolute constant $c$. Obviously, $\{a^k\}$ and $\{b^k\}$ are weakly equicontinuous in $C([0,T]; l^2)$, and hence are relatively compact in $C([0,T]; l^2_{\mathrm w})$ by the Ascoli-Arzela theorem. Therefore, there exists a subsequence $\{\left(a^{k_n}, b^{k_n}\right)\}$ and a weakly continuous pair of functions $(a(t), b(t))$ in $l^2\times l^2$ such that
\begin{equation}\notag
a^{k_n} \to a, \ \ b^{k_n} \to b \ \ \ \mbox{as} \ \ \ k_n\to \infty \ \ \mbox{in} \ \ C([0,T]; l^2_{\mathrm w}).
\end{equation}
Hence,
\begin{equation}\notag
a_j^{k_n}(t) \to a_j(t), \ \ b_j^{k_n}(t) \to b_j(t) \ \ \ \mbox{as} \ \ \ k_n\to \infty \ \ \mbox{for all}\  j\geq 0 \ \ \mbox{and} \ \ t\in[0,T].
\end{equation}
In particular, we have $a(0)=a^0$ and $b(0)=b^0$. 

It is left to show that the limit $(a,b)$ is a solution of (\ref{sys-2}). By (\ref{sys-truncated}), we have
\begin{equation}\notag
\begin{split}
a_j^{k_n}(t)=&\ a_j^{k_n}(0)+\int_0^t  \left(-\lambda_j^{\theta}a_j^{k_n}(\tau)a_{j+1}^{k_n}(\tau)+\lambda_{j-1}^{\theta}(a_{j-1}^{k_n}(\tau))^2\right) \, d\tau\\
&+\int_0^t \left(-\lambda_{j}^{\theta}b_j^{k_n}(\tau)b_{j+1}^{k_n}(\tau)+\lambda_{j-1}^{\theta}(b_{j-1}^{k_n}(\tau))^2+f_j\right) \, d\tau,\\
b_j^{k_n}(t)=&\ b_j^{k_n}(0)+\int_0^t  \left(\lambda_j^{\theta}a_j^{k_n}(\tau)b_{j+1}^{k_n}(\tau)-\lambda_{j}^{\theta}b_{j}^{k_n}(\tau)a_{j+1}^{k_n}(\tau)\right) \, d\tau,
\end{split}
\end{equation}
for $0\leq j\leq k_n-1$. After taking the limit $k_n\to\infty$ we get
\begin{equation}\notag
\begin{split}
a_j(t)=&\ a_j(0)+\int_0^t  \left(-\lambda_j^{\theta}a_j(\tau)a_{j+1}(\tau)+\lambda_{j-1}^{\theta}(a_{j-1}(\tau))^2\right) \, d\tau\\
&+\int_0^t \left(-\lambda_{j}^{\theta}b_j(\tau)b_{j+1}(\tau)+\lambda_{j-1}^{\theta}(b_{j-1}(\tau))^2+f_j\right) \, d\tau,\\
b_j(t)=&\ b_j(0)+\int_0^t  \left(\lambda_j^{\theta}a_j(\tau)b_{j+1}(\tau)-\lambda_{j}^{\theta}b_{j}(\tau)a_{j+1}(\tau)\right) \, d\tau.
\end{split}
\end{equation}
Furthermore, since $a_j$ and $b_j$ are continuous for all $j\geq 0$, the integral form implies that $a_j\in C^1[0,T]$ and $b_j\in C^1[0,T]$. Thus $(a, b)$ solves (\ref{sys-2}). 

\cbdu

\bigskip

\section{Local well-posedness}
\label{sec-well}

System (\ref{sys-2}) and (\ref{sys-1}) are locally well-posed in $H^s$ with $s\geq \theta$. Namely, we prove:

\begin{Theorem}\label{thm-strong}
Let $s\geq \theta$. Assume $(a^0, b^0) \in H^s\times H^s$ and $f\in H^s$. There exists a time $T>0$ such that system (\ref{sys-2}) with initial data $(a^0, b^0)$ has a unique solution $(a(t), b(t))$ in $H^s\times H^s$ on $[0,T]$.
\end{Theorem}

\pf
The local existence follows from a standard argument and a priori estimate for solutions in $H^s$. We only show the a priori estimate here. Multiplying the $a_j$ equation of (\ref{sys-2}) by $\lambda_j^{2s}a_j$, the $b_j$ equation by $\lambda_j^{2s}b_j$, and taking the sum over $j\geq 0$, we obtain
\begin{equation}\notag
\begin{split}
\frac{1}{2}\frac{d}{dt}\sum_{j=0}^\infty \left(\lambda_j^{2s}a_j^2+\lambda_j^{2s} b_j^2\right)
=&(\lambda^{2s}-1)\sum_{j=0}^\infty \lambda_j^{\theta+2s} a_j^2a_{j+1}\\
& +(\lambda^{2s}-1)\sum_{j=0}^\infty \lambda_j^{\theta+2s} b_j^2a_{j+1}
 +\sum_{j=0}^\infty \lambda_j^{2s} a_jf_j.
\end{split}
\end{equation}
Applying Cauchy-Schwartz's inequality gives
\begin{equation}\notag
\begin{split}
\sum_{j=0}^\infty \lambda_j^{\theta+2s} a_j^2a_{j+1}=&\ \lambda^{-\frac13(\theta+2s)}\sum_{j=0}^\infty \lambda_j^{\frac23(\theta+2s)} a_j^2 \cdot \lambda_{j+1}^{\frac13(\theta+2s)}a_{j+1}\\
\leq&\ \frac23\lambda^{-\frac13(\theta+2s)}\sum_{j=0}^\infty \lambda_j^{\theta+2s} a_j^3
+\frac13\lambda^{-\frac13(\theta+2s)}\sum_{j=0}^\infty \lambda_{j+1}^{\theta+2s} a_{j+1}^3\\
\leq &\ \lambda^{-\frac13(\theta+2s)}\left(\sum_{j=0}^\infty \lambda_j^{\frac23(\theta+2s)} a_j^2\right)^{\frac32},
\end{split}
\end{equation}
and similarly
\begin{equation}\notag
\begin{split}
\sum_{j=0}^\infty \lambda_j^{\theta+2s} b_j^2a_{j+1}=&\ \lambda^{-\frac13(\theta+2s)}\sum_{j=0}^\infty \lambda_j^{\frac23(\theta+2s)} b_j^2 \cdot \lambda_{j+1}^{\frac13(\theta+2s)}a_{j+1}\\
\leq&\ \frac23\lambda^{-\frac13(\theta+2s)}\sum_{j=0}^\infty \lambda_j^{\theta+2s} b_j^3
+\frac13\lambda^{-\frac13(\theta+2s)}\sum_{j=0}^\infty \lambda_{j+1}^{\theta+2s} a_{j+1}^3\\
\leq &\ \frac23\lambda^{-\frac13(\theta+2s)}\left(\sum_{j=0}^\infty \lambda_j^{\frac23(\theta+2s)} b_j^2\right)^{\frac32}\\
&+\frac13\lambda^{-\frac13(\theta+2s)}\left(\sum_{j=0}^\infty \lambda_j^{\frac23(\theta+2s)} a_j^2\right)^{\frac32}.
\end{split}
\end{equation}
We also have
\begin{equation}\notag
\sum_{j=0}^\infty \lambda_j^{2s} a_jf_j\leq \sum_{j=0}^\infty \lambda_j^{2s} a_j^2+ \sum_{j=0}^\infty \lambda_j^{2s}f_j^2.
\end{equation}
Therefore, putting the last four estimates together yields
\begin{equation}\notag
\begin{split}
\frac{1}{2}\frac{d}{dt}\sum_{j=0}^\infty \left(\lambda_j^{2s}a_j^2+\lambda_j^{2s} b_j^2\right)
\leq &\frac43(\lambda^{2s}-1)  \lambda^{-\frac13(\theta+2s)}\left(\sum_{j=0}^\infty \lambda_j^{\frac23(\theta+2s)} a_j^2\right)^{\frac32}\\
&+\frac23(\lambda^{2s}-1)  \lambda^{-\frac13(\theta+2s)}\left(\sum_{j=0}^\infty \lambda_j^{\frac23(\theta+2s)} b_j^2\right)^{\frac32}\\
&+\sum_{j=0}^\infty \lambda_j^{2s} a_j^2+ \sum_{j=0}^\infty \lambda_j^{2s}f_j^2.
\end{split}
\end{equation}
Since $s\geq \theta$, we have $\frac23(\theta+2s)\leq 2s$ and hence
\begin{equation}\notag
\begin{split}
\frac{1}{2}\frac{d}{dt}\sum_{j=0}^\infty \left(\lambda_j^{2s}a_j^2+\lambda_j^{2s} b_j^2\right)
\leq &C(\lambda, s, \theta)\left(\sum_{j=0}^\infty( \lambda_j^{2s} a_j^2+\lambda_j^{2s} b_j^2)\right)^{\frac32}\\
&+\sum_{j=0}^\infty \lambda_j^{2s} a_j^2+ \sum_{j=0}^\infty \lambda_j^{2s}f_j^2
\end{split}
\end{equation}
for a constant $C(\lambda, s, \theta)$. 
Thus, for $f\in H^s$, there exists a time $T>0$ and another constant $C(\lambda, s, \theta, f)$ such that 
\begin{equation}\notag
\|a(t)\|_{H^s}^2+\|b(t)\|_{H^s}^2 \leq C(\lambda, s, \theta, f)\left(\|a(0)\|_{H^s}^2+\|b(0)\|_{H^s}^2\right), \ \ \ t\in[0,T].
\end{equation}
The uniqueness can be obtained by the classical approach, i.e. establishing a Gr\"owall's inequality for the difference of two solutions starting from the same data. We omit the details here.

\cbdu

\begin{Theorem}\label{thm-strong2}
Let $s\geq \theta$. Assume $(a^0, b^0) \in H^s\times H^s$ and $f\in H^s$. There exists a time $T>0$ such that system (\ref{sys-1}) with initial data $(a^0, b^0)$ has a unique solution $(a(t), b(t))$ in $H^s\times H^s$ on $[0,T]$.
\end{Theorem}

The proof follows similarly as that of Theorem \ref{thm-strong}.

\begin{Remark}
Recall $\theta=\frac{5-\delta}{2}$ and $\delta\in[0,3]$. 
When $\theta=\frac52$ ($\delta=0$) the local well-posedness of (\ref{sys-2}) and (\ref{sys-1}) holds in $H^s$ with $s\geq \frac52$ which is consistent with the well-posedness result for the ideal MHD system. On the other hand, if $\theta=1$ ($\delta=3$), (\ref{sys-2}) and (\ref{sys-1}) are locally well-posed in $H^s$ with $s\geq 1$. 
\end{Remark}

\bigskip

\section{Finite time blow-up}
\label{sec-blow}

In this section we show that a solution of the dyadic ideal MHD model (\ref{sys-2}) starting from positive initial data develops blow-up at a finite time in the space $H^s$ with $s>\frac13\theta$. 
For dyadic models, we are aware of two approaches of proving finite time blow-up. For the forced dyadic Euler system the authors of \cite{CFP2} made use of  properties of the fixed point and provided a dynamical approach to show the development of blow-up. Another method is to construct a suitable Lyapunov function for the underlying system that blows up at finite time, for instance, see \cite{Dai-21, JL}. In particular, the model studied in \cite{JL} contains both the DN type and Obukov type of nonlinearities, and exhibits both forward and backward energy cascade mechanisms. Consequently, the solutions considered there are not known to be positive. However, the authors demonstrated blow-up by developing a special Lyapunov function approach. Their approach has certain flexibility and is robust.
 In the context of model (\ref{sys-2}), the interactions of the fluid and magnetic field make it harder to explore the energy transfer from shell to shell; as a result, a dynamical approach to show blow-up is out of reach. Therefore we choose to take the path of a Lyapunov function. Nevertheless, it seems quite challenging to construct finite time blow-up solution for system (\ref{sys-1}) by applying the aforementioned approaches. In fact, the mechanism of both forward and backward energy cascades may prevent finite time blow-up.

The result of blow-up for (\ref{sys-2}) is stated below.

\begin{Theorem}\label{thm-blow-2}
Let $\theta>0$ and $f_0\geq 0$.
The solution $(a(t), b(t))$ of (\ref{sys-2}) with positive initial data develops blow-up at a finite time in the $H^s$ norm with $s>\frac13\theta$. 
\end{Theorem}
\pf
The plan is to show that the quantity 
\[E_s(t):= \|a(t)\|_{H^s}^2+\|b(t)\|_{H^s}^2 \]
with $s>\frac13\theta$ is not locally integrable, which will be achieved through a contradiction argument.
 
We introduce the rescaled variables $w_j=\lambda_j^\theta a_j$ and $z_j=\lambda_j^\theta b_j$. It follows from (\ref{sys-2}) that $w=(w_0, w_1, w_2, ...)$ and $z=(z_0, z_1, z_2, ...)$ satisfy 
\begin{equation}\label{sys-wz}
\begin{split}
w_j'=&\ -\lambda^{-\theta} w_jw_{j+1}+\lambda^\theta w_{j-1}^2-\lambda^{-\theta} z_jz_{j+1}+\lambda^\theta z_{j-1}^2+\lambda_j^\theta f_j,\\
z_j'=&\ \lambda^{-\theta}w_j z_{j+1}-\lambda^{-\theta} z_jw_{j+1},
\end{split}
\end{equation}
for $j\geq 0$, with $w_{-1}=z_{-1}=0$, $f_0>0$ and $f_j=0$ for $j\geq 1$.
Consider the quantities 
\begin{equation}\notag
\begin{split}
\phi(t)=&\ \sum_{j=0}^\infty \left(\lambda_j^{-\gamma} w_j^2(t)+\lambda_j^{-\gamma} z_j^2(t)\right),\\
\psi(t)=&\ \sum_{j=0}^\infty \left(\lambda_j^{-\gamma} w_j(t)+c_0\lambda_j^{-\gamma} z_j(t)\right),
\end{split}
\end{equation}
for appropriate constants $\gamma>0$ and $0<c_0<1$ to be specified later.

Assume $E_s(t)$
is locally integrable for $s>\frac13\theta$.  Notice that 
\[0<\phi(t)=\sum_{j=0}^\infty \left(\lambda_j^{2\theta-\gamma} a_j^2(t)+\lambda_j^{2\theta-\gamma} b_j^2(t)\right)\leq E_s(t) \]
provided
\begin{equation}\label{para-condition-1}
2\theta-\gamma\leq 2s,
\end{equation}
and hence $\phi(t)$ is locally integrable.
We will proceed to show that $\psi(t)$ is locally integrable as well and in the same time it  satisfies a Riccati type of inequality for an appropriate value of $\gamma$, which leads to an obvious contradiction.  

First of all, by Cauchy-Schwarz's inequality, we have
\begin{equation}\label{ineq-phi-psi}
\begin{split}
\psi^2(t)=&\ \left(\sum_{j=0}^\infty \left(\lambda_j^{-\gamma} w_j(t)+c_0\lambda_j^{-\gamma} z_j(t)\right)\right)^2\\
\leq &\  \left(\sum_{j=0}^\infty \lambda_j^{-\gamma}\right)\sum_{j=0}^\infty \lambda_j^{-\gamma} \left(w_j(t)+ c_0z_j(t)\right)^2\\
\leq &\ \frac{2}{1-\lambda^{-\gamma}}\sum_{j=0}^\infty \lambda_j^{-\gamma} \left(w_j^2(t)+c_0^2 z_j^2(t)\right)\\
\leq &\ \frac{2c_0^2\phi(t)}{1-\lambda^{-\gamma}}.
\end{split}
\end{equation}
As an immediate consequence, $\psi(t)$ is locally integrable. 

On the other hand, straightforward computation based on (\ref{sys-wz}) shows 
\begin{equation}\label{energy-wz}
\begin{split}
\frac{d}{dt}\psi(t)=& -\sum_{j=0}^\infty\lambda^{-\theta}\lambda_j^{-\gamma} w_jw_{j+1}+\sum_{j=0}^\infty\lambda^{\theta}\lambda_j^{-\gamma} w_{j-1}^2\\
&-\sum_{j=0}^\infty\lambda^{-\theta}\lambda_j^{-\gamma} z_jz_{j+1}
+\sum_{j=0}^\infty\lambda^{\theta}\lambda_j^{-\gamma} z_{j-1}^2\\
&+c_0\sum_{j=0}^\infty\lambda^{-\theta}\lambda_j^{-\gamma} w_jz_{j+1}
-c_0\sum_{j=0}^\infty\lambda^{-\theta}\lambda_j^{-\gamma} z_jw_{j+1}+f_0.
\end{split}
\end{equation}
Since $\phi(t)$ is locally integrable, it is finite almost everywhere. Notice that the infinite sums on the right hand side of (\ref{energy-wz}) are in the order of $\phi(t)$  and hence are defined almost everywhere. 
It is clear that
\begin{equation}\label{est-positive-1}
\sum_{j=0}^\infty\lambda^{\theta}\lambda_j^{-\gamma} w_{j-1}^2+\sum_{j=0}^\infty\lambda^{\theta}\lambda_j^{-\gamma} z_{j-1}^2=\lambda^{\theta-\gamma} \phi(t).
\end{equation}
Applying Cauchy-Schwarz's inequality, we estimate the infinite sums as follows
\begin{equation}\notag
\begin{split}
\sum_{j=0}^\infty\lambda^{-\theta}\lambda_j^{-\gamma} w_jw_{j+1}
=&\ \lambda^{\frac{\gamma}{2}-\theta} \sum_{j=0}^\infty\left(\lambda_j^{-\frac{\gamma}{2}} w_j\right)\left(\lambda_{j+1}^{-\frac{\gamma}{2}}w_{j+1}\right)\\
\leq &\ \lambda^{\frac{\gamma}{2}-\theta} \left(\sum_{j=0}^\infty\lambda_j^{-\gamma} w_j^2\right)^{\frac12} \left(\sum_{j=0}^\infty\lambda_{j+1}^{-\gamma} w_{j+1}^2\right)^{\frac12},\\
\end{split}
\end{equation}
and similarly
\begin{equation}\notag
\begin{split}
\sum_{j=0}^\infty\lambda^{-\theta}\lambda_j^{-\gamma} z_jz_{j+1}\leq&\ \lambda^{\frac{\gamma}{2}-\theta} \left(\sum_{j=0}^\infty\lambda_j^{-\gamma} z_j^2\right)^{\frac12} \left(\sum_{j=0}^\infty\lambda_{j+1}^{-\gamma} z_{j+1}^2\right)^{\frac12}, \\
\sum_{j=0}^\infty\lambda^{-\theta}\lambda_j^{-\gamma} z_jw_{j+1}\leq&\ \lambda^{\frac{\gamma}{2}-\theta} \left(\sum_{j=0}^\infty\lambda_j^{-\gamma} z_j^2\right)^{\frac12} \left(\sum_{j=0}^\infty\lambda_{j+1}^{-\gamma} w_{j+1}^2\right)^{\frac12},
\end{split}
\end{equation}
\begin{equation}\notag
-\sum_{j=0}^\infty\lambda^{-\theta}\lambda_j^{-\gamma} w_jz_{j+1}\leq \lambda^{\frac{\gamma}{2}-\theta} \left(\sum_{j=0}^\infty\lambda_j^{-\gamma} w_j^2\right)^{\frac12} \left(\sum_{j=0}^\infty\lambda_{j+1}^{-\gamma} z_{j+1}^2\right)^{\frac12}.
\end{equation}
Combining the last four inequalities yields
\begin{equation}\label{est-positive-4}
\begin{split}
&\sum_{j=0}^\infty\lambda^{-\theta}\lambda_j^{-\gamma} w_jw_{j+1}+\sum_{j=0}^\infty\lambda^{-\theta}\lambda_j^{-\gamma} z_jz_{j+1}\\
&-c_0\sum_{j=0}^\infty\lambda^{-\theta}\lambda_j^{-\gamma} w_jz_{j+1}+c_0\sum_{j=0}^\infty\lambda^{-\theta}\lambda_j^{-\gamma} z_jw_{j+1}\\
\leq &\ (1+2c_0) \lambda^{\frac{\gamma}{2}-\theta} \phi(t).
\end{split}
\end{equation}
Since $a_j(t)>0$ and $b_j(t)>0$ for all $j\geq 0$ and $t\geq 0$, it is also true that 
$w_j(t)>0$ and $z_j(t)>0$ for all $j\geq 0$ and $t\geq 0$. Therefore, we deduce from (\ref{energy-wz}), (\ref{est-positive-1}), (\ref{est-positive-4}) and (\ref{ineq-phi-psi}) that
\begin{equation}\notag
\begin{split}
\frac{d}{dt}\psi(t)
\geq &\ \left(\lambda^{\theta-\gamma}-(1+2c_0)\lambda^{\frac{\gamma}{2}-\theta}\right) \phi(t) +f_0\\
\geq &\ \frac{1-\lambda^{-\gamma}}{2c_0^2}\left(\lambda^{\theta-\gamma}-(1+2c_0)\lambda^{\frac{\gamma}{2}-\theta}\right) \psi^2(t)+f_0
\end{split}
\end{equation}
provided $\lambda^{\theta-\gamma}-(1+2c_0)\lambda^{\frac{\gamma}{2}-\theta}>0$. 
It follows from the Riccati type of inequality that $\psi(t)$ becomes infinity at a finite time if 
\begin{equation}\label{cond-initial}
\frac{1-\lambda^{-\gamma}}{2c_0^2}\left(\lambda^{\theta-\gamma}-(1+2c_0)\lambda^{\frac{\gamma}{2}-\theta}\right) \psi^2(0)+f_0>0.
\end{equation}
Therefore, the contradiction is achieved. 

In the end, we summarize all the conditions on the parameters
\begin{equation}\notag
2\theta-\gamma\leq 2s, \ \ \lambda^{\theta-\gamma}-(1+2c_0)\lambda^{\frac{\gamma}{2}-\theta}>0.
\end{equation}
In order to have the latter one satisfied, we can choose 
\[\theta-\gamma>\frac{\gamma}{2}-\theta, \ \ 0<c_0<\frac12\lambda^{2\theta-\frac32\gamma}-\frac12.\]
Thus, $\gamma$ can be chosen as 
\[2\theta-2s\leq \gamma<\frac43\theta.\]
Notice that $[2\theta-2s, \frac43\theta)$ is not empty given the assumption $s>\frac13\theta$. We also note condition (\ref{cond-initial}) is satisfied for any positive initial data and forcing $f_0\geq0$. 
\cbdu

\medskip

We note the importance of the critical threshold index $\frac13\theta$ in terms of the energy equality. For the original ideal MHD (\ref{mhda})-(\ref{mhdc}), it was proved in \cite{CKS} and \cite{KL} that solutions conserving energy need to have the minimum regularity of $\frac13$, which is analogous to Onsager's critical regularity index for the Euler equation. 
We note that $\frac13\theta=\frac56$ for $\theta=\frac52$ (and $\delta=0$); and the space $H^{\frac56}$ has the same scaling as $B^{\frac13}_{3, \infty}$ and $B^{\frac13}_{3, c(\mathbb N)}$ in three dimensions, which have Onsager's critical scaling proved in \cite{CKS} and \cite{KL}.

\bigskip

\section{The set of fixed points}
\label{sec-fixed}

We consider (\ref{sys-2}) and (\ref{sys-1}) with the forcing $f=(f_0, f_1, f_2, ...)$ satisfying 
\[ f_0>0; \ \ f_j=0, \ \ j\geq 1. \]
We will show the existence of fixed points of both systems (\ref{sys-2}) and (\ref{sys-1}); moreover, we provide an explicit form of them.

In order to simplify notations, we will rescale the stationary system of (\ref{sys-2}) and (\ref{sys-1}). Define 
\[A_j=\lambda^{-\frac16\theta}f_0^{-\frac12}\lambda_j^{\frac13\theta}a_j, \ \ \ B_j=\lambda^{-\frac16\theta}f_0^{-\frac12}\lambda_j^{\frac13\theta}b_j.\]
Let $\frac{d}{dt}a_j=\frac{d}{dt}b_j=0$ for all $j\geq 0$.  System (\ref{sys-2}) reduces to
\begin{subequations}
\begin{align}
A_{j-1}^2-A_jA_{j+1}+B_{j-1}^2-B_jB_{j+1}=&\ 0, \ \ j\geq 1, \label{stat-2-mhd1}\\
A_jB_{j+1}-B_jA_{j+1} =&\ 0, \ \ j\geq 1, \label{stat-2-mhd2}\\
-A_0A_1-B_0B_1=& -1, \label{stat-2-mhd3}\\
A_0B_1-B_0A_1=&\ 0. \label{stat-2-mhd4}
\end{align}
\end{subequations}
\begin{Lemma}\label{le-ex-fixed}
The solutions $\{(A_j,B_j)\}$ of (\ref{stat-2-mhd1})-(\ref{stat-2-mhd4}) satisfy
\[A_j=A_0, \ \ B_j=B_0, \ \ \forall \ \ j\geq1,\]
for some constants $A_0$ and $B_0$ with the constraint
\[A_0^2+B_0^2=1.\]
\end{Lemma}
\pf
The two equations $(\ref{stat-2-mhd2})$ and $(\ref{stat-2-mhd4})$ imply 
\begin{equation}\label{stat-21-2}
\frac{A_{j+1}}{A_j}=\frac{B_{j+1}}{B_j}=c_j, \ \ \forall j\geq0,
\end{equation}
for certain constant $c_j$. It follows that $A_{j+1}=c_jA_j$ and $B_{j+1}=c_jB_j$ for all $j\geq 0$. Then $(\ref{stat-2-mhd1})$ can be rewritten as
\begin{equation}\label{stat-22-2}
\left(A_j^2+B_j^2\right)\left(c_{j-1}^{-2}-c_j\right)=0, \ \ \forall j\geq 1.
\end{equation}
It follows
\begin{equation}\notag
c_j=c_{j-1}^{-2}, \ \ \forall j\geq1.
\end{equation}
As a consequence of the recursive relation, we have
\begin{equation}\label{stat-24-2}
c_j= c_0^{(-2)^j}, \ \ \ j\geq0.
\end{equation}
Therefore, (\ref{stat-21-2}) together with (\ref{stat-24-2}) gives rise to 
\begin{equation}\notag
\begin{split}
A_{j+1}=&\ c_jA_j=\cdot\cdot\cdot =c_jc_{j-1}\cdot\cdot\cdot c_1c_0A_0\\
=&\ c_0^{(-2)^j+(-2)^{j-1}+\cdot\cdot\cdot +(-2)+1}A_0\\
=&\ c_0^{\left(1-(-2)^{j+1}\right)/3}A_0,
\end{split}
\end{equation}
and 
\begin{equation}\notag
B_{j+1}= c_0^{\left(1-(-2)^{j+1}\right)/3}B_0.
\end{equation}
We further simplify the previous form and obtain
\begin{equation}\label{stat-25-2}
A_{j+1}=
\begin{cases}
c_0^{\frac{1+2^{j+1}}{3}}A_0, \ \ j \ \mbox{is even,}\\
c_0^{\frac{1-2^{j+1}}{3}}A_0, \ \ j \ \mbox{is odd,}
\end{cases}
B_{j+1}=
\begin{cases}
c_0^{\frac{1+2^{j+1}}{3}}B_0, \ \ j \ \mbox{is even,}\\
c_0^{\frac{1-2^{j+1}}{3}}B_0, \ \ j \ \mbox{is odd.}
\end{cases}
\end{equation}
Consequently, we claim that $c_0$ has to be $1$, since the sequences $\{a_j\}$ and $\{b_j\}$ are bounded. Thus, (\ref{stat-24-2}) implies that $c_j=1$ for all $j\geq 1$. Hence, $A_j=A_0$ and $B_j=B_0$ for all $j\geq 1$. In addition, it follows from $(\ref{stat-2-mhd3})$ that 
\begin{equation}\notag
A_0^2+B_0^2=1.
\end{equation}
It concludes the proof of the lemma.

\cbdu

\medskip

With the same rescaling as above, the stationary system of (\ref{sys-1}) can be written as
\begin{subequations}
\begin{align}
A_{j-1}^2-A_jA_{j+1}-B_{j-1}^2+B_jB_{j+1}=&\ 0, \ \ j\geq 1, \label{stat-1-mhd1}\\
-A_jB_{j+1}+B_jA_{j+1} =&\ 0, \ \ j\geq 1, \label{stat-1-mhd2}\\
-A_0A_1+B_0B_1=& -1, \label{stat-1-mhd3}\\
-A_0B_1+B_0A_1=&\ 0. \label{stat-1-mhd4}
\end{align}
\end{subequations}

\begin{Lemma}\label{le-ex-fixed1}
The solutions $\{(A_j,B_j)\}$ of (\ref{stat-1-mhd1})-(\ref{stat-1-mhd4}) satisfy
\[A_j=A_0, \ \ B_j=B_0, \ \ \forall \ \ j\geq1,\]
for some constants $A_0$ and $B_0$ with the constraint
\[A_0^2-B_0^2=1.\]
\end{Lemma}
\pf
The proof is analogous as that of Lemma \ref{le-ex-fixed}. 
It follows from equations $(\ref{stat-1-mhd2})$ and $(\ref{stat-1-mhd4})$ that 
\[A_{j+1}=c_jA_j, \ \ \ B_{j+1}=c_jB_j\] 
for certain constant $c_j$ and all $j\geq 0$. 
Thus we obtain from $(\ref{stat-1-mhd1})$ that
\begin{equation}\label{stat-22-2}
\left(A_j^2-B_j^2\right)\left(c_{j-1}^{-2}-c_j\right)=0, \ \ \forall j\geq 1.
\end{equation}
Therefore,
\begin{equation}\notag
c_j=c_{j-1}^{-2}, \ \ \mbox{or} \ \ A_j^2-B_j^2=0 \ \ \forall j\geq1.
\end{equation}
Thanks to (\ref{stat-1-mhd2}) and (\ref{stat-1-mhd4}), we observe if $A_j^2-B_j^2=0$ for one $j$, it has to be true for all $j\geq 0$. However,  (\ref{stat-1-mhd3}) implies $c_0(A_0^2-B_0^2)=1$. Consequently, $A_0^2-B_0^2\neq 0$ and hence $A_j^2-B_j^2\neq 0$ for all $j\geq 0$. We conclude the only possibility is 
$c_j=c_{j-1}^{-2}$ for all $j\geq0$. 
It then follows from the same argument from the proof of Lemma \ref{le-ex-fixed} that $c_j=1$ for all $j\geq 1$. Thus we have $A_j=A_0$ and $B_j=B_0$ for all $j\geq 1$. In the end, $(\ref{stat-1-mhd3})$ implies that 
\begin{equation}\notag
A_0^2-B_0^2=1.
\end{equation}

\cbdu

\bigskip

\section{Linear instability}
\label{sec-stab}

This section concerns the long time behavior of solutions, specifically, whether a solution of (\ref{sys-2}) converges to a fixed point or not. Evidently, from Lemma \ref{le-ex-fixed} and Lemma \ref{le-ex-fixed1}, we know that a fixed point of (\ref{sys-2}) or (\ref{sys-1}) is not unique. 
Without loss of generality, we choose $f_0=\lambda^{-\frac13\theta}$ such that $\lambda^{\frac16\theta} f_0^{\frac12}=1$; thus, the steady state $(\bar a, \bar b)$ can be written as 
\begin{equation}\label{exact-fix}
\bar a_j=A_0\lambda_j^{-\frac13\theta}, \ \ 
\bar b_j=B_0\lambda_j^{-\frac13\theta}, \ \ \forall \ \ j\geq 0.
\end{equation}
Consider the perturbation of the steady state,
\begin{equation}\label{def-diff-omega}
\begin{split}
a_j(t)=&\ A_0\lambda_j^{-\frac13\theta}+\epsilon \omega_j(t), \ \ \ j\geq0,\\
b_j(t)=&\ B_0\lambda_j^{-\frac13\theta}+\epsilon \zeta_j(t), \ \ \ j\geq0.
\end{split}
\end{equation}
Based on (\ref{sys-2}), we derive the equations for $\omega_j$ and $\zeta_j$
\begin{equation}\label{eq-omega}
\begin{split}
\omega_j'= &\ A_0\lambda_j^{\frac23\theta}\left(2\lambda^{-\frac23\theta} \omega_{j-1}- \lambda^{-\frac13\theta} \omega_j-\omega_{j+1}\right)\\
&+B_0\lambda_j^{\frac23\theta}\left(2\lambda^{-\frac23\theta} \zeta_{j-1}- \lambda^{-\frac13\theta} \zeta_j-\zeta_{j+1}\right)\\
&- \epsilon\lambda_j^\theta \omega_j\omega_{j+1}+\epsilon\lambda_{j-1}^\theta \omega_{j-1}^2
- \epsilon\lambda_j^\theta \zeta_j\zeta_{j+1}+\epsilon\lambda_{j-1}^\theta \zeta_{j-1}^2, \\ 
\zeta_j'= &\ B_0\lambda_j^{\frac23\theta} \left(\lambda^{-\frac13\theta} \omega_j-\omega_{j+1}\right)
-A_0\lambda_j^{\frac23\theta} \left(\lambda^{-\frac13\theta} \zeta_j-\zeta_{j+1}\right)\\
&+ \epsilon\lambda_j^\theta \omega_j\zeta_{j+1}-\epsilon \lambda_j^\theta \zeta_j\omega_{j+1}
\end{split}
\end{equation}
for $j\geq 0$ and $\omega_{-1}=\zeta_{-1}=0$.
We focus on the two special cases: 
\begin{itemize}
\item [(i)] $A_0=1, B_0=0$: the fixed point corresponds to the fixed point of the dyadic Euler equation. We show linear instability with stable velocity component and unstable magnetic field. 
\item [(ii)] $A_0=-1, B_0=0$: the fixed point also corresponds to a steady state of the dyadic Euler equation. We show linear instability with both unstable velocity and unstable magnetic field. 
\end{itemize}

The main results described in (i) and (ii) are presented below.
\begin{Theorem}\label{thm-stab1}
Let $(\bar a, \bar b)$ be the fixed point of (\ref{sys-2}) with $A_0=1$ and $B_0=0$ (hence $\bar b=0$). For the linearized system about $(\bar a, \bar b)$, there are no positive eigenvalues corresponding to the velocity component; while every real number is an eigenvalue for the magnetic field linearized equation. 
\end{Theorem}
\pf
When $A_0=1$ and $B_0=0$, it follows from (\ref{eq-omega}) that the linearized system is
\begin{subequations}
\begin{align}
\omega_j'= &\ \lambda_j^{\frac23\theta}\left(2\lambda^{-\frac23\theta} \omega_{j-1}- \lambda^{-\frac13\theta} \omega_j-\omega_{j+1}\right) \label{linear-omega1}\\
\zeta_j'= &
-\lambda_j^{\frac23\theta} \left(\lambda^{-\frac13\theta} \zeta_j-\zeta_{j+1}\right) \label{linear-omega2}
\end{align}
\end{subequations}
for $j\geq 0$ and $\omega_{-1}=\zeta_{-1}=0$. Notice that $\omega_j$ and $\zeta_j$ are decoupled in the linearized system, and (\ref{linear-omega1}) is the linearized equation of dyadic Euler model about the fixed point $\bar a_j=\lambda_j^{-\frac13\theta}$. We will apply the continued fraction approach from \cite{CFP1} to (\ref{linear-omega1})-(\ref{linear-omega2}). 
We look for a solution to (\ref{linear-omega1}) and (\ref{linear-omega2}) in the form
\[\omega_j=c_j e^{p t}, \ \ \zeta_j=d_je^{q t}\]
with real values of $p$ and $q$.  
Inserting the form $\omega_j=c_j e^{p t}$ to (\ref{linear-omega1}) yields
\begin{equation}\label{eq-c}
c_{j+1}+\left(p\lambda_j^{-\frac23\theta}+\lambda^{-\frac13\theta}\right)c_j-2\lambda^{-\frac23\theta}c_{j-1}=0
\end{equation}
for $j\geq 0$, with $c_{-1}=0$. Inserting $\zeta_j=d_je^{q t}$ to (\ref{linear-omega2}) gives
\begin{equation}\label{eq-d}
d_{j+1}=\left(q\lambda_j^{-\frac23\theta}+\lambda^{-\frac13\theta}\right)d_j, \ \ j\geq 0
\end{equation}
with $d_{-1}=0$. 

Following the arguments of \cite{CFP1}, one can deduce from the iterative equation (\ref{eq-c}) that it is impossible to have $p\geq0$. Therefore, for a real and negative $p$, we can choose $c_0=1$, and hence 
\[c_{j}=\left(p\lambda_{j-1}^{-\frac23\theta}+\lambda^{-\frac13\theta}\right)\left(p\lambda_{j-2}^{-\frac23\theta}+\lambda^{-\frac13\theta}\right)\cdot\cdot\cdot\left(p\lambda_{1}^{-\frac23\theta}+\lambda^{-\frac13\theta}\right)\lambda^{-\frac13\theta}.\]
It implies that 
\[c_j\sim \lambda_j^{-\frac13\theta}, \ \ \forall \ \ j\geq J \ \ \mbox{for some} \ \ J>0.\]
Therefore, the sequence $\{c_j\}$ has finite $H^s$ norm for $s<\frac13\theta$, and in particular has finite $l^2$ norm.

On the other hand, it follows from (\ref{eq-d}) that
\begin{equation}\notag
d_{j+1}=\alpha_j\alpha_{j-1}\alpha_0d_0, \ \ \mbox{with} \ \ \alpha_j=q\lambda_j^{-\frac23\theta}+\lambda^{-\frac13\theta}.
\end{equation}
As a consequence, we know 
\begin{equation}\label{limit}
\lim_{j\to\infty}\alpha_j=\lambda^{-\frac13\theta}, \ \ \ \ \lim_{j\to\infty}d_j=0,
\end{equation}
which hold for any $q\in \mathbb R$.
In view of (\ref{limit}), we also have 
\begin{equation}\notag
d_j\sim \lambda_j^{-\frac13\theta}, \ \ \ \mbox{for large} \ \ j.
\end{equation}
Thus the sequence $\{c_j\}$ has finite $H^s$ norm for $s<\frac13\theta$. Hence the $H^s$ norm of $\{\zeta_j\}$ is finite for all the time if $q\leq 0$ and grows exponentially if $q>0$.

\cbdu

\medskip

\begin{Theorem}\label{thm-instab1}
Let $(\bar a, \bar b)$ be the fixed point of (\ref{sys-2}) with $A_0=-1$ and $B_0=0$. There are no negative eigenvalue to the velocity equation in the linearized system about such $(\bar a, \bar b)$. Every real number is an eigenvalue for the magnetic field equation in the linearized system. Thus the linearized system is unstable.
\end{Theorem}
\pf
Taking $A_0=-1$ and $B_0=0$ in (\ref{eq-omega}) gives the linearized system 
\begin{equation}\label{linear2}
\begin{split}
\omega_j'= &- \lambda_j^{\frac23\theta}\left(2\lambda^{-\frac23\theta} \omega_{j-1}- \lambda^{-\frac13\theta} \omega_j-\omega_{j+1}\right) \\
\zeta_j'= &\ \lambda_j^{\frac23\theta} \left(\lambda^{-\frac13\theta} \zeta_j-\zeta_{j+1}\right) 
\end{split}
\end{equation}
Similarly, looking for solutions of the linearized system in the form
\[\omega_j=c_j e^{p t}, \ \ \zeta_j=d_je^{q t},\]
we have
\begin{equation}\label{eq-c2}
c_{j+1}+\left(\lambda^{-\frac13\theta}-p\lambda_j^{-\frac23\theta}\right)c_j-2\lambda^{-\frac23\theta}c_{j-1}=0,
\end{equation}
\begin{equation}\label{eq-d2}
d_{j+1}=\left(\lambda^{-\frac13\theta}-q\lambda_j^{-\frac23\theta}\right)d_j,
\end{equation}
for $j\geq 0$. Now denote $\alpha_j=\lambda^{-\frac13\theta}-p\lambda_j^{-\frac23\theta}$. Following the analysis of \cite{CFP1}, we infer from (\ref{eq-c2}) that for $j\geq 1$
\begin{equation}\label{eq-continue}
p \lambda_j^{-\frac23\theta}-\lambda^{-\frac13\theta}=-\alpha_j=[\alpha_{j+1}, \alpha_{j+2}, ...]
\end{equation}
with the continued fraction defined as
\begin{equation}\notag
[\alpha_{j+1}, \alpha_{j+2}, ...]=\frac{1}{\lambda^{\frac32\theta}\alpha_{j+1}+\frac{\lambda^{\frac32\theta}}{\lambda^{\frac32\theta}\alpha_{j+2}+\frac{\lambda^{\frac32\theta}}{\cdot\cdot\cdot}}}
\end{equation}
We note that $\alpha_j>0$ for all $j\geq J$ for some large enough $J$, and hence the right hand side of (\ref{eq-continue}) is positive. That implies $p$ needs to be positive. Therefore, we know that the velocity component is unstable about the steady state with $A_0=-1$ and $B_0=0$. 

It follows from (\ref{eq-d2}) that
\begin{equation}\label{limit}
\lim_{j\to\infty}\alpha_j=\lambda^{-\frac13\theta}, \ \ \ \ \lim_{j\to\infty}d_j=0,
\end{equation}
for any $q\in \mathbb R$. Thus the magnetic field is unstable due to a similar analysis as in the proof of Theorem \ref{thm-stab1}.

\cbdu

\begin{Theorem}\label{thm-instab3}
Let $(\bar a, \bar b)$ be the fixed point of (\ref{sys-1}) with $A_0=-1$ and $B_0=0$. There are no negative eigenvalue to the velocity equation in the linearized system about such $(\bar a, \bar b)$. Every real number is an eigenvalue for the magnetic field equation in the linearized system. Thus the linearized system is unstable.
\end{Theorem}
\pf
In analogy with previous analysis for system (\ref{sys-2}), considering perturbation 
\begin{equation}\notag
a_j(t)= A_0\lambda_j^{-\frac13\theta}+\epsilon \omega_j(t), \ \ \ 
b_j(t)= B_0\lambda_j^{-\frac13\theta}+\epsilon \zeta_j(t), \ \ \ j\geq0
\end{equation}
for system (\ref{sys-1}), the pair $(\omega_j, \zeta_j)$ satisfies 
\begin{equation}\notag
\begin{split}
\omega_j'= &\ A_0\lambda_j^{\frac23\theta}\left(2\lambda^{-\frac23\theta} \omega_{j-1}- \lambda^{-\frac13\theta} \omega_j-\omega_{j+1}\right)\\
&-B_0\lambda_j^{\frac23\theta}\left(2\lambda^{-\frac23\theta} \zeta_{j-1}- \lambda^{-\frac13\theta} \zeta_j-\zeta_{j+1}\right)\\
&- \epsilon\lambda_j^\theta \omega_j\omega_{j+1}+\epsilon\lambda_{j-1}^\theta \omega_{j-1}^2
+ \epsilon\lambda_j^\theta \zeta_j\zeta_{j+1}-\epsilon\lambda_{j-1}^\theta \zeta_{j-1}^2, \\ 
\zeta_j'= &-B_0\lambda_j^{\frac23\theta} \left(\lambda^{-\frac13\theta} \omega_j-\omega_{j+1}\right)
+A_0\lambda_j^{\frac23\theta} \left(\lambda^{-\frac13\theta} \zeta_j-\zeta_{j+1}\right)\\
&- \epsilon\lambda_j^\theta \omega_j\zeta_{j+1}+\epsilon \lambda_j^\theta \zeta_j\omega_{j+1}
\end{split}
\end{equation}
for $j\geq 0$ and $\omega_{-1}=\zeta_{-1}=0$. If $A_0=1$ and $B_0=0$ the linearized system with $\epsilon=0$ is 
\begin{equation}\label{linear3}
\begin{split}
\omega_j'= &\ \lambda_j^{\frac23\theta}\left(2\lambda^{-\frac23\theta} \omega_{j-1}- \lambda^{-\frac13\theta} \omega_j-\omega_{j+1}\right),\\
\zeta_j'= &\ \lambda_j^{\frac23\theta} \left(\lambda^{-\frac13\theta} \zeta_j-\zeta_{j+1}\right),
\end{split}
\end{equation}
while if $A_0=-1$ and $B_0=0$ the linearized system is
\begin{equation}\label{linear4}
\begin{split}
\omega_j'= & - \lambda_j^{\frac23\theta}\left(2\lambda^{-\frac23\theta} \omega_{j-1}- \lambda^{-\frac13\theta} \omega_j-\omega_{j+1}\right),\\
\zeta_j'= & - \lambda_j^{\frac23\theta} \left(\lambda^{-\frac13\theta} \zeta_j-\zeta_{j+1}\right).
\end{split}
\end{equation}
In view of the form of the linearized systems (\ref{linear-omega1})-(\ref{linear-omega2}), (\ref{linear2}), (\ref{linear3}) and (\ref{linear4}), we note that the instability results stated in Theorem \ref{thm-stab1} and Thorem \ref{thm-instab1} hold as well for the linearized systems of (\ref{sys-1}) about the fixed point with $A_0=1$ and $B_0=0$ and the fixed point with $A_0=-1$ and $B_0=0$ respectively.

\cbdu

\bigskip

\section*{Acknowledgement}
The work of M. Dai is partially supported by NSF Grants DMS--1815069 and DMS--2009422; the work of S. Friedlander is partially supported by NSF Grant DMS--1613135. S. Friedlander is grateful to IAS for its hospitality in 2020-2021.

\bigskip

\section*{Conflict of interest statement}
The authors certify that there is no conflict of interest.

\end{document}